\documentclass{amsart}

\usepackage{amscd,amsthm,xy,amssymb}
\xyoption{all}
\newdir{ >}{{}*!/-5pt/@{>}}

\newcommand{\Ab}{\mathfrak{Ab}}

\newcommand{\C}{\mathfrak{C}}

\newcommand{\fin}{\mathfrak{Fin}}

\newcommand{\F}{\mathcal{F}}

\newcommand{\cha}{\rm{Ch}^{\ge 0}}
\newcommand{\iso}{\overset{\cong}{\to}}
\newcommand{\weq}{\overset{\sim}{\longrightarrow}}
\newcommand{\cofi}{\rightarrowtail}
\newcommand{\fib}{\twoheadrightarrow}
\newcommand{\trfib}{\overset{\sim}{\twoheadrightarrow}}
\newcommand{\pe}{\hat{p}}

\newcommand{\bp}{\tilde{P}}
\newcommand{\bq}{\tilde{Q}}

\DeclareMathOperator{\ho}{Ho}
\DeclareMathOperator{\sets}{Sets}

\DeclareMathOperator{\coker}{coker}

\DeclareMathOperator{\ass}{Rings}

\DeclareMathOperator{\comm}{Comm}
\DeclareMathOperator{\sur}{Sur}
\DeclareMathOperator{\map}{Map}
\DeclareMathOperator{\sign}{sign}

\newcommand{\Z}{\mathbb{Z}}
\newcommand{\Q}{\mathbb{Q}}

\theoremstyle{plain}

\newtheorem{them}{Theorem}[section]
\newtheorem{coro}[them]{Corollary}
\newtheorem{lemm}[them]{Lemma}
\newtheorem{propo}[them]{Proposition}

\theoremstyle{definition}

\newtheorem{defi}[them]{Definition}
\newtheorem{rema}[them]{Remark}

\newtheorem{nota}[them]{Notation}
\newtheorem{powers}[them]{Simplicial powers and cosimplicial homotopies}
\newtheorem{product}[them]{Product structure}

\begin{document}

\title[Dold-Kan correspondence for rings]{Cosimplicial versus DG-rings: a version of the Dold-Kan correspondence}
\author[J.L. Castiglioni]{Jos\'e Luis Castiglioni*}
\address{Departamento de Matem\'atica\\
         Facultad de Cs. Exactas\\
           Calle 50 y 115\\
         (1900) La Plata\\
            Argentina.}

\email{jlc@mate.unlp.edu.ar}

\author[G. Corti\~nas]{Guillermo Corti\~nas**}

\address{Departamento de Matem\'atica\\
Ciudad Universitaria\\
 Pabell\'on 1\\
(1428) Buenos Aires\\
 Argentina.}

\email{gcorti@dm.uba.ar}

\begin{abstract}
The (dual) Dold-Kan correspondence says that there is an equivalence of
categories $K:\cha\to \Ab^\Delta$ between nonnegatively graded cochain complexes and cosimplicial abelian groups, which is inverse to the normalization 
functor. We show that the restriction of $K$ to $DG$-rings can be equipped
with an associative product and that the resulting functor 
$DGR^*\to\ass^\Delta$, although not itself an equivalence, does induce one
at the level of homotopy categories. In other words both $DGR^*$ and $\ass^\Delta$ are Quillen closed model categories and the total left derived functor of $K$
is an equivalence:
\begin{equation*}
\mathbb{L}K:\ho DGR^*\weq\ho \ass^\Delta
\end{equation*}
The dual of this result for chain $DG$ and simplicial rings was obtained independently by 
S. Schwede and B. Shipley through different methods ({\it Equivalences of 
monoidal model categories}.
Algebraic and Geometric Topology 3 (2003), 287-334). Our proof is based
on a functor $Q:DGR^*\to \ass^\Delta$, naturally homotopy equivalent to $K$, and which
preserves the closed model structure. It also has other interesting
applications. For example, we use $Q$ to prove a noncommutative version of the 
Hochschild-Kostant-Rosenberg and Loday-Quillen theorems. Our version applies to the cyclic
module $[n] \mapsto \coprod^n_R S$ that arises from a homomorphism $R\to S$ of not necessarily commutative rings, using the coproduct $\coprod_R$ of associative $R$-algebras.
As another application of the properties of $Q$, we obtain a simple, braid-free 
description of a product on the tensor power $S^{\otimes_R^n}$ originally 
defined by P. Nuss using braids ({\it Noncommutative descent and nonabelian cohomology,} 
K-theory {\bf 12} (1997) 23-74.). 
\end{abstract}

\thanks{(*), (**): Both authors were partially supported by grant
UBACyT X066. The second author is a CONICET researcher and an ICTP
associate}

\maketitle
\section{Introduction}
The (dual) Dold-Kan correspondence is an equivalence between the category
$\cha$ of nonnegatively graded cochain complexes of abelian groups and the
category $\Ab^\Delta$ of cosimplicial abelian groups. This equivalence is
defined by a pair of inverse functors
\begin{equation}
N:\Ab^{\Delta}\leftrightarrows \cha:K
\end{equation}
Here $N$ is the normalized or Moore complex (see \eqref{normac} below). 
The functor $K$ is described
in \cite{we}, 8.4.4; if $A=(A,d)\in\cha$ and $n\ge 0$, then 
\begin{equation}
K^nA=\bigoplus_{i=0}^n\binom{n}{i}A^i\cong\bigoplus_{i=0}^n A^i\otimes\Lambda^i\Z^n
\end{equation}
If in addition $A$ happens to be a $DG$-ring, then $K^nA$ can be equipped
with a product, namely that coming from the tensor product of rings 
$A\otimes\Lambda \Z^n$:
\begin{equation}\label{intrinus}
(a\otimes x)(b\otimes y)=ab\otimes x\land y.
\end{equation}
This product actually makes 
$[n]\mapsto K^nA$ into a cosimplicial ring (see 5.3). 
Thus $K$ can be viewed as a functor from $DG$- to cosimplicial rings:
\begin{equation}
K:DGR^*\to \ass^\Delta,\qquad A\mapsto KA
\end{equation}
Note that for all $n$, $K^nA$ is a nilpotent extension of $A^0$.
As there are cosimplicial rings which are not codimensionwise 
nilpotent extensions of constant cosimplicial rings, $A\mapsto KA$ is not
a category equivalence. However we prove (Theorem \ref{comain}) that it induces
one upon inverting weak equivalences. Precisely, $K$
carries quasi-isomorphisms to maps inducing an isomorphism at the cohomotopy
level, and therefore induces a functor $\mathbb{L}K$
between the localizations $\ho DGR^*$ and $\ho\ass^\Delta$ obtained by 
formally inverting such maps, and we prove that $\mathbb{L}K$ is an equivalence:
\begin{equation}\label{equintro}
\mathbb{L}K:\ho DGR^*\weq\ho\ass^\Delta
\end{equation} 
The dual of this result, that is, the equivalence between the homotopy categories
of chain $DG$ and simplicial rings, was obtained independently by Schwede and 
Shipley through different methods (see \cite{ship} and also Remark \ref{LneQ}
below). 

To prove \eqref{equintro} we use Quillen's formalism of closed model categories
\cite{qui}. We consider in each of $DGR^*$ and $\ass^\Delta$ a closed model 
structure,
in which weak equivalences are as above, fibrations are surjective maps
and cofibrations are appropiately defined to fit Quillen's axioms.
There is a technical problem in that the functor $K$ does not preserve
cofibrations. To get around this, we replace $K$ by a certain functor $Q$.
As is the case of the Dold-Kan functor, $Q$ too is defined for all cochain 
complexes $A$, even if they may not be $DG$-rings. If $A\in\cha$ then
\begin{equation}
Q^nA=\bigoplus_{i=0}^\infty A^i\otimes T^i(\Z^n)
\end{equation}
We show that any set map $\alpha:[n]\to [m]$ induces a group homomorphism
$Q^nA\to Q^mA$, so that $[n]\mapsto Q^nA$ is not only a functor on $\Delta$
but on the larger category $\fin$ with the same objects, where a 
homomorphism $[n]\to [m]$ is just
any set map. The projection $T\Z^n\to \Lambda\Z^n$ induces a homomorphism
\begin{equation}
\pe:QA\weq KA
\end{equation}
We show $\pe$ induces an isomorphism of cohomotopy groups.
If moreover $A$ is a $DG$-ring, $Q^nA$ has an obvious product coming from
$A\otimes T\Z^n$; however this product is not well-behaved with respect
to the $\fin$ nor the cosimplicial structure. In order to get a $\fin$-ring
we perturb the product by a Hochschild $2$-cocycle 
$f:A^*\otimes T^*V\to A^{*+1}\otimes T^{*+1}V$. We obtain a
product $\circ$ of the form 
\begin{equation}\label{circo}
(a\otimes x)\circ (b\otimes y)=ab\otimes xy+f(a\otimes x,b\otimes y)
\end{equation}
For a definition of $f$ see \eqref{fedo} below. It turns out that the
map $\pe$ is a ring homomorphism (see 5.3). 
This implies that the derived functors of $K$ and of the 
functor $\tilde{Q}$ obtained from $Q$ by restriction of its $\fin$-structure to
a cosimplicial one, are isomorphic (see \ref{qpreserve}):
\begin{equation}
\mathbb{L}\tilde{Q}\cong\mathbb{L}K.
\end{equation}
We show further that $\mathbb{L}\tilde{Q}$ is an equivalence. We deduce this
from the stronger result (Theorem \ref{main1}) that $\tilde{Q}$ is the
left adjoint of a Quillen equivalence 
(as defined in Hovey's book \cite{hov}, 1.3.12).

Next we review other results obtained in this paper.
As mentioned above, for $A\in\cha$, $QA$ is not only a cosimplicial
group but a $\fin$-group. In particular the cyclic permutation
$t_n:=(0\dots n):[n]\to [n]$ acts on $Q^nA$, and we may view $QA$ as a cyclic
module in the sense of \cite{we}, 9.6.1. Consider the associated 
normalized mixed complex $(NQA,\mu,B)$. We show that there is a weak
equivalence of mixed complexes
\begin{equation}\label{loqu}
(A,0,d)\weq (NQA,\mu,B)
\end{equation}
In particular these two mixed complexes have the same Hochschild homology:
\begin{equation}\label{isoshu}
A^*\cong H_*(NQA,\mu)
\end{equation}
If $A$ happens to be a $DG$-ring then the shuffle product induces a graded
ring structure on $H_*(NQA,\mu)$; we show in \ref{propi} that \eqref{isoshu}
is a ring isomorphism for the product of $A$ and the shuffle product of
$H_*(NQA,\mu)$.

A specially interesting case is that of the $DG$-ring of noncommutative
differential forms $\Omega_RS$ relative to a ring homomorphism $R\to S$
(as defined in \cite{cq}). We show in \ref{qomega} that $Q\Omega_RS$ is the
coproduct $\fin$-ring:
\begin{equation}\label{introqomega}
Q\Omega_RS\cong{\coprod}_RS:[n]\mapsto{\coprod_{i=0}^n}{}_{{}_R}S
\end{equation} 
In particular, by \eqref{isoshu}, there is an isomorphism of graded rings 
\begin{equation}\label{introhkr}
\Omega_RS\iso H_*(N{\coprod}_RS,\mu)
\end{equation}
The particular case of \eqref{introhkr} when $R$ is commutative and $R\to S$ is
central and flat was proved in 1994 by Guccione, Guccione and Majadas \cite{ggm}.
More generally, by \eqref{loqu} we have a mixed complex equivalence
\begin{equation}\label{introlq}
(\Omega_RS,0,d)\weq (N{\coprod}_RS,\mu,B)
\end{equation} 
We view \eqref{introhkr} and \eqref{introlq} as noncommutative versions
of the Hochschild-Kostant-Rosenberg and Loday-Quillen theorems \cite{we}
9.4.13, 9.8.7.

As another application, we give a simple formulation
for a product structure defined by Nuss \cite{nuss} on each term
of the Amitsur complex associated to a homomorphism $R\to S$ of not necessarily
commutative rings $R$ and $S$:
\begin{equation}
{\bigotimes}_RS:[n]\mapsto{\bigotimes_{i=0}^n}{}_{{}_R}S
\end{equation}
Nuss constructs his product using tools from the theory of quantum groups. 
We show
here (see Section \ref{conu}) that the 
canonical Dold-Kan isomorphism maps the product \eqref{intrinus} to that
defined by Nuss. Thus 
\begin{equation}
K\Omega_RS=KN({\bigotimes}_RS)\cong{\bigotimes}_RS
\end{equation}
is an isomorphism of cosimplicial rings.

The remainder of this paper is organized as follows. Basic notations are fixed
in Section \ref{linear}. In Section \ref{q} the functor $Q$ is defined.
The homotopy equivalence of the cosimplicial groups $KA$ and $QA$ as
well as that of the mixed complexes \eqref{loqu} is proved in Section
\ref{compaqk}. In Section \ref{Qismonoidal} we show that the functor
$Q:\cha\to \Ab^\fin$ is strong monoidal (\ref{Qisstrong}). We use this to
introduce, for $A\in DGR^*$, the product \eqref{circo} on $QA$ 
(5.3). 
The graded ring isomorphism
\eqref{isoshu} is proved in Section \ref{comshu}. The isomorphism
\eqref{introqomega} and its corollaries \eqref{introhkr} and \eqref{introlq} 
are
proved in Section \ref{nchkr}. The reformulation of Nuss' product is
the subject of Section \ref{conu}. 
In Section \ref{equintro} we prove that $\tilde{Q}$ is the left adjoint of a
Quillen equivalence (Theorem \ref{main1}) and deduce from this that 
$\mathbb{L}K$ is a category equivalence (Corollary \ref{comain}).

\section{Cochain complexes and cosimplicial abelian groups}
\label{linear}
We write $\Delta$ for the simplicial category, and $\fin$ for the category with
the same objects as $\Delta$, but where the homomorphisms $[n]\to [m]$ are
just the set maps. The inclusion
\begin{equation*}
\hom_\Delta([n],[m])\subset\map([n],[m])=\hom_\fin([n],[m])
\end{equation*}
gives a faithful embedding $\Delta\subset \fin$.
If $I$ and $\C$ are categories, we shall write $\C^I$ to
denote the category of functors $I\to \C$, to which we refer as $I$-objects
of $\C$. If $C:I\to \C$ is an $I$-object, we write $C^i$ for $C(i)$. We use
the same letter for a map $\alpha:[n]\to [m]\in I$ as for its image under $C$. The canonical embedding $\Delta\subset\fin$ mentioned above
makes $[n]\mapsto [n]$ into a cosimplicial object of $\fin$. We write
$\partial_i:[n]\to [n+1]$, $i=0,\dots, n+1$ and $\mu_j:[n]\to [n-1]$,
$j=0,\dots,n-1$, for the coface and codegeneracy maps. We also
consider the map $\mu_n:[n]\to [n-1]$ defined by
\begin{equation}\label{muen}
\mu_n(i)=\left\{\begin{matrix}i&\text{ if }&i<n\\
                            0&\text{ if }&i=n\end{matrix} \right .
\end{equation}
One checks that $d_i:=\mu_i:[n]\to [n-1]$, $i=0,\dots,n$ and
$s_j=\partial_{j+1}:[n]\to [n+1]$, $j=0,\dots,n$ satisfy the simplicial
identities, with the $d_i$ as faces and the $s_i$ as degeneracies.
Thus there is a functor $\Delta^{op}\to \fin$, $[n]\mapsto [n]$.
Moreover the cyclic
permutation $t_n=(0\dots n):[n]\to [n]$ extends this simplicial
structure to a cyclic one (see \cite{we}, 9.6.3). Composing with these functors
and with the inclusion $\Delta\subset\fin$ mentioned above
we have a canonical way of regarding any $\fin$-object in a
category $\C$ as either a cosimplicial, a simplicial, or a cyclic
object.

If $\C$ is a category with finite
coproducts, and $A\in\C$, we write $\coprod A$ for the functor

\begin{equation}\label{copra}
\coprod A:\fin\to\C,\qquad [n]\mapsto \coprod_{i=0}^nA
\end{equation}
Here $\coprod$ may be replaced by whatever sign denotes the coproduct of
$\C$; for example if $\C$ is abelian, we write $\oplus A$ for $\coprod A$.

If $A=\oplus_{n=0}^\infty A_n$ and $B=\oplus_{n=0}^\infty B_n$ are graded abelian
groups, we write
\begin{equation}\label{boxtimes}
A\boxtimes B:=\oplus_{n=0}^\infty A_n\otimes B_n
\end{equation}
If $A,B$ are graded $I$-abelian groups, we put $A\boxtimes B$
for the graded $I$-abelian group $i\mapsto A^i\boxtimes B^i$.

\section{The functor $Q$}\label{q}
We are going to define a functor $Q:Ch^{\ge 0}\to \Ab^\fin$; first we
need some auxiliary constructions.
Write $V:=\ker(\oplus \Z\to \Z)$ for the kernel of the canonical map to the
constant $\fin$-abelian group, and $\{e_i:0\le i\le n\}$ for the canonical
basis of $\oplus_{i=0}^n\Z$. Put $v_i=e_i-e_0$, $0\le i\le n$. Note $v_0=0$
and $\{v_1,\dots,v_n\}$ is a basis of $V^n$. The action of a map $\alpha:[n]\to [m]\in\fin$ on $V$ is given by

\begin{equation}\label{alphav}
\alpha v_i=v_{\alpha (i)}-v_{\alpha(0)}\qquad (0\le i\le n)
\end{equation}

Applying to $V$ the tensor algebra functor
$T$ in each codimension yields a graded $\fin$-ring $TV$. If $A=(A,d)\in\cha$,
we put
\begin{equation}
Q^nA:=A\boxtimes TV^n
\end{equation}
If $\alpha:[n]\to [m]\in \fin$, we set
\begin{equation}\label{alphaQ}
\alpha(a\otimes x)=a\otimes\alpha x+da\otimes v_{\alpha(0)}\alpha x
\end{equation}
If $\beta:[m]\to[p]\in \fin$, then
\begin{align*}
\beta(\alpha (a\otimes x))=&a\otimes\beta\alpha x+da\otimes v_{\beta(0)}
\beta\alpha x+
da\otimes\beta(v_{\alpha (0)})\beta\alpha x\\
                          =&(\beta\alpha)(a\otimes x)
\end{align*}
Thus $QA$ is a $\fin$-abelian group, and
$Q:\cha\to\Ab^{\fin}$ a functor.
We have a filtration on $QA$ by
$\fin$-subgroups, given by
\begin{equation}\label{filq}
\F_nQA=\bigoplus_{i=n}^\infty A^i\otimes T^iV
\end{equation}
The associated graded $\fin$-abelian group is $G_{\F}QA=A\boxtimes TV$.

\section{Comparison between $Q$ and the Dold-Kan functor $K$}
\label{compaqk}
The Dold-Kan correspondence is a pair of inverse functors (see \cite{we} 8.4):
\begin{equation*}
K:\cha\leftrightarrows\Ab^{\Delta}:N
\end{equation*}
If $C\in\Ab^{\Delta}$
then $NC$ can be
equivalently described as the normalized complex or as the Moore complex:
\begin{equation}\label{normac}
N^nC=C^n/\sum_{i=1}^n\partial_iC^{n-1}\cong\cap_{i=0}^{n-1}\ker(\mu_i:C^n\to C^{n-1})
\end{equation}
In either version the coboundary map $N^nC\to N^{n+1}C$ is induced by
\begin{equation}\label{partial}
\partial=\sum_{i=0}^n(-1)^i\partial_i.
\end{equation}
In the first version this is
the same map as that induced by $\partial_0$.
A description of the inverse functor $K$ (in the simplicial case)
is given in \cite{we}, 8.4.4, and another in \cite{kar}, 1.5.
Here is yet another.
Let $\Lambda V$ be the exterior algebra, $p:TV\to \Lambda V$
the canonical projection. One checks that $\ker(1\otimes p)\subset QA$
is a $\fin$-subgroup. Thus
\begin{equation}\label{ka}
K^*A:=A\boxtimes\Lambda V^*
\end{equation}
inherits a $\fin$-structure. Moreover
\begin{equation}\label{mape}
\pe:=1\otimes p:QA\fib KA
\end{equation}
is a natural
surjection of $\fin$-abelian groups. To see that the resulting cosimplicial
abelian group $KA$ is indeed the same as (i.e. is naturally isomorphic to) that of \cite{we}, it suffices to show that $NKA=A$. Put
$$V_j^n=\bigoplus_{i\ne j}\Z v_i\subset \bigoplus_{i=1}^n\Z v_i=V^n$$
We have
\begin{align*}
NK^nA=&A\boxtimes\Lambda V^n/\sum_{i=1}^nA\boxtimes\partial_i(\Lambda V^n)\\
     =&A\boxtimes(\Lambda V^n/\sum_{i=1}^n\Lambda(V_i^n))\\
     =&A^n\otimes v_1\land\dots\land v_n\cong A^n
\end{align*}
Furthermore it is clear that the coboundary map induced by $\partial_0$
is $d:A^*\to A^{*+1}$. Thus our $KA$ is the same cosimplicial abelian
group as that of \cite{we}. But since in our construction $KA$ has a
$\fin$-structure, we may also regard it as a simplicial or cyclic abelian
group. From our definition of faces and degeneracies, it is clear that
the normalized complex of $KA$ considered as a simplicial group has the
abelian group $N^nKA=A^n$ in each dimension. One checks that the alternating 
sum $\mu$ of
the faces induces the trivial boundary. Thus the normalized chain
complex of the simplicial group $KA$ is $(A,0)$.
Consider the Connes operator $B:NQ^*A\to NQ^{*+1}A$,
\begin{equation}\label{B}
B=\partial_0\circ\sum_{i=0}^n(-1)^{ni}t_n^i
\end{equation}
We show in \ref{kequivq} below that $\pe B=D\pe$, where $D:=(n+1)d$ on
$A^n$. Hence we have a map of mixed complexes
\begin{equation}\label{pemix}
\pe:(NQA,\mu,B)\to (A,0,D).
\end{equation}
We shall see in \ref{kequivq} below that \eqref{pemix} is a {\it rational
equivalence of mixed complexes}. We recall that a map of mixed complexes
is an equivalence if it induces an isomorphism at the level of Hochschild
homology; this automatically implies it also induces an isomorphism at the
level of cyclic, periodic cyclic and negative cyclic homologies.
In \ref{kequivq} we also consider the map
\begin{equation}\label{elele}
l:A\to NQA,\qquad l(a)=a\otimes\sum_{\sigma\in S_n}\sign(\sigma)
v_{\sigma 1}\dots v_{\sigma n}
\end{equation}
We show in Theorem \ref{kequivq} below that $l$ is an integral equivalence

\[l:(A,0,d)\weq (NQA,\mu,B)\]

\smallskip

\begin{rema}\label{rescape}
Note that if $A$ is a complex of $\Q$-vectorspaces, then $\pe$
can be rescaled as $(1/n!)\pe$ on $NQ^nA$ to give a mixed complex
map $(NQA,\mu,B)\to (A,0,d)$ which is left inverse to $l$.
\end{rema}

\smallskip

\begin{them}\label{kequivq}
Let $A$ be a cochain complex of abelian groups, $\pe:QA\fib KA$
the map of $\fin$-abelian groups defined in \eqref{mape} above. Then:

\noindent{i)} There are a natural cochain map $j:(A,d)\to (NQA,\partial)$ such that $\pe j=1_{A}$ and a natural cochain homotopy $h:N^*QA\to N^{*-1}QA$
such that $[h,\partial]=1-j\pe$.

\noindent{ii)} The map \eqref{pemix} is a rational equivalence of mixed
complexes. On the other hand the map \eqref{elele} is a natural integral
equivalence $l:(A,0,d)\to (NQA,\mu,B)$.
\end{them}
\begin{proof}
First we compute $NQA$. A similar argument as that given in Section
\ref{compaqk} to
compute $NKA$, shows that
\begin{equation}\label{nqdef}
N^nQA=A\boxtimes (TV^n/\sum_{i=1}^nTV^n_i).
\end{equation}
On the other hand we have a canonical identification between
the $r$th tensor power of $V^n=\Z^n$ and the free abelian group on the set
of all maps $\{1,\dots,r\}\to\{1,\dots,n\}$:
\begin{equation}\label{identensorfun}
T^rV^n\cong \Z[\map(\{1,\dots,r\},\{1,\dots,n\})]
\end{equation}
Using \eqref{identensorfun}, $T^rV^n/\sum_{j=1}^n T^rV^n_j$ becomes
the free module on all surjective maps $\{1,\dots,r\}\to\{1,\dots,n\}$;
we get
\begin{equation}\label{nqa}
N^nQA=A\boxtimes \Z[\sur_{*,n}]=\bigoplus_{r=n}^\infty A^r\otimes
\Z[\sur_{r,n}]
\end{equation}
Here $\sur_{p,q}$ is the set of all surjections $\{1, ..., p\} \to \{1, ..., q\}$. 
Note that in particular $\sur_{n,n} = S_n$, the symmetric group on $n$ letters.
To prove i), regard $NQA$ as a cochain complex.
We may view $NQA$ as the direct sum total complex of a second quadrant
double complex
\begin{equation*}
C^{p,q} = 
\left\{
\begin{array}{cc}
A & \textrm{ if $p = q = 0$,}\\
A^q\otimes \Z[\sur_{q,q+p}] & \textrm{ if $(p,q) \neq (0,0)$.}
\end{array}
\right.
\end{equation*}
Here $1\otimes \partial_0$ and $d\otimes v_1\partial_0$ are respectively
the horizontal and the vertical coboundary operators. The filtration
\eqref{filq} is the row filtration.
If we regard $A=NKA$ as a double cochain complex concentrated in the zero column, then $\pe$ becomes a map of double complexes.
By definition, $\pe=1\otimes p$; at the $n$th row, $p$ is a map:
\begin{equation}\label{gepi}
p:\Z[\sur_{n,*}]\fib \Z[n].
\end{equation}
The only nonzero component of $p$ is $p(\sigma)=\sign(\sigma)$.
We claim \eqref{gepi} is a cochain homotopy equivalence.
To prove this note first that because both $\Z[\sur_{n,*}]$ and $\Z[n]$ are
complexes of free abelian
groups, to show $p$ is a homotopy equivalence it suffices to check it is
a quasi-isomorphism. Next note that
\begin{align}\label{cohotv}
H^*(\Z[\sur_{n,*}])=&H^*(NT^nV)=\pi^*(T^nV)\\
=&T^n\pi^*(V)\nonumber\\
=&T^nH^*(NV)\nonumber\\
=&T^nH^*(\Z[1])=\Z[n].\nonumber
\end{align}
Thus, to prove $p$ is a cochain equivalence
it suffices to show that
\begin{equation}\label{kerp}
\ker(p:\Z[S_n]\to \Z)=\partial_0(\Z[\sur_{n,n-1}])
\end{equation}
The inclusion $\supset$ of \eqref{kerp} holds
because $p$ is a cochain map. To prove the
other inclusion, proceed as follows. First note the identification
\begin{equation*}
\Z[S_n]\cong \bigoplus_{\sigma\in S_n}\Z v_{\sigma 1}\dots v_{\sigma n}
\end{equation*}
Next observe that the kernel of $p$ is generated by elements of the form
\[
\dots v_1\dots v_i\dots+\dots v_i \dots v_1\dots
 \equiv -\partial_0(\dots v_{i-1}\dots v_{i-1}\dots)\qquad (i>1)
\]
Here congruence is taken modulo $\sum_{j\ge 1}\partial_jTV$.
Thus $p$ is a surjective homotopy equivalence, as claimed.
Therefore we may choose a cochain map
$j':\Z[n]\to \Z[\sur_{n,*}]$ such that $pj'=1$ and a cochain homotopy
$h': \Z[\sur_{n,*}]\to \Z[\sur_{n,*-1}]$ such that $[h',\partial_0]=1-pj'$.
One checks that the following maps satisfy the requirements of part i) of
the theorem:
\begin{align*}
j:=&1\otimes j'+(1\otimes h')((1\otimes j')d-d\otimes v_1\partial_0j')\\
h:=&(1\otimes h'-(1\otimes h')(d\otimes v_1\partial_0)(1\otimes h'))(1\otimes j'p-1)
\end{align*}
Next we prove part ii). Observe the face maps of $NQA$ are of the form
$1\otimes\mu_i$ where $\mu_i$ is the face map in $TV$.
Hence we have a direct sum decomposition of chain complexes
\begin{equation}\label{nqmu}
(NQA,\mu)=\bigoplus_{n=0}^\infty A^n\otimes(\Z[\sur_{n,*}],\mu)
\end{equation}
The homology version of the argument used in \eqref{cohotv} shows
that $$H_*(\Z[\sur_{n,*}])=\Z[n].$$
In particular $L_n:=\ker(\mu:\Z[S_n]\to\Z[\sur_{n,n-1}])$ is free of rank one.
By definition, to prove $\pe$ is a rational mixed complex equivalence, we must
prove that $\pe\mu=0$, which is straighforward, that $\pe B=D\pe$, which we leave for
later, and finally that $\pe=1\otimes p:(NQA,\mu)\to (A,0)$ is a rational chain equivalence, which in turn reduces to proving $p(L_n)\ne 0$ for $n\ge 1$.
Consider the element
\begin{equation}\label{symm}
\epsilon_n:=\sum_{\sigma\in S_n}\sign(\sigma)\sigma\in\Z[S_n]
\end{equation}
We have $p(\epsilon_n)=n!$; one checks further that $\epsilon_n\in L_n$.
It follows that $\pe:(NQA,\mu)\to (A,0)$ is a rational equivalence, as
we had to prove. Moreover, as every coefficient of $\epsilon_n$ is
invertible, and $L_n$ has rank one, we have $L_n=\Z\epsilon_n$. It follows that
the map $l':\Z[n]\to (\Z[\sur_{n,*}],\mu)$ which sends $1\in \Z$ to $\epsilon_n$ is a quasi-isomorphism, whence a homotopy equivalence. To finish the proof,
we must show that $ld=Bl$ and $\pe B=D\pe$. Both of these follow once one has
proven the formula \eqref{formuB} below, which in turn is derived from the
identities \eqref{lemB}, which are proved by induction. The inclusion $\{1\} \subseteq \{1,..., n+1\}$ together with the map $\{1,\dots,n\}\to\{1,\dots,n+1\}$, $i \mapsto i+1$, define a 
bijection $\{1\}\coprod\{1,\dots,n\} \to \{1,\dots,n+1\}$. We identify $\{1\}\coprod\{1,\dots,n\}=\{1,\dots,n+1\}$ using this bijection.
If $\sigma\in S_n$, we denote by
$1\coprod\sigma$ the coproduct map.
\begin{equation}\label{formuB}
B(a\otimes \sigma)=da\otimes\sum_{i=0}^n(-1)^{in}(1\dots n+1)^i(1\coprod\sigma)
\end{equation}

\begin{align}\label{lemB}
t_n^i(v_j)=&\begin{cases}v_{i+j}-v_i&\text{ if }i\le n-j\nonumber\\
                        v_{p-1}-v_i& \text{if }i=n-j+p\quad j\ge p\ge 1
\end{cases}\\
t_n^i(a\otimes x)=&a\otimes t_n^ix+da\otimes v_{i}t_n^ix\\
B(a\otimes x)=&da\otimes\sum_{i=0}^n(-1)^{in}v_{i+1}\partial_0t^ix\nonumber
\end{align}
\end{proof}
\begin{nota} Let $B=(B,d)\in\cha$. Put $PB^n=B^n\oplus B^{n-1}\oplus B^n$.
Equip $PB$ with the coboundary operator $\partial:PB^*\to PB^{*+1}$
given by the matrix
\begin{equation*}
\partial=\left[\begin{matrix}d&0&0\\
                             1&-d&-1\\
                             0&0&d\end{matrix}\right]
\end{equation*}
We note $PB$ comes equipped with a natural map
$\epsilon=(\epsilon_0,\epsilon_1):PB\to B\oplus B$, and that two maps
$f_0,f_1:A\to B$
are cochain homotopic if and only if there exists a cochain homomorphism
$H:A\to PB$ such that
$\epsilon H=(f_0,f_1)$.\qed
\end{nota}

\medskip	

The next corollary says that, for $A,B\in\cha$, every cosimplicial
map $f:QA\to QB$ has a canonically associated cochain map $\bar{f}$,
such that $NQ\bar{f}$ and $Nf$ are naturally homotopic. Moreover if
$f=Qg$, then $\bar{f}=g$.

\medskip

\begin{coro}\label{barf}
Let $A,B\in\cha$. Consider the functors
\begin{align*}
({\cha})^{op}\times\cha &\to\Ab\\
(A,B)\mapsto & \hom_{\Ab^\Delta}(QA,QB)\\
(A,B)\mapsto & \hom_{\cha}(NQA,PNQB).\\
\end{align*}
There are two natural transformations
\begin{align*}
\bar{\ \ :}\hom_{\Ab^\Delta}(QA,QB)\to&\hom_{\cha}(A,B),\\
H:\hom_{\Ab^\Delta}(QA,QB)\to & \hom_{\cha}(NQA,PNQB).\\
\end{align*}
These are such that $\overline{Qg}=g$ and
that the following diagram commutes
\begin{align}
\xymatrix{&\hom_{\cha}(NQA,PNQB)\ar[d]^(.4){\epsilon}\\
          \hom_{\Ab^\Delta}(QA,QB)\ar[ur]^H\ar[r]&\hom_{\cha}(NQA,NQB\oplus NQB)}& \\
\nonumber \xymatrix{f \ar@{|->}[rrrrr]& & & & & (Nf, NQ \bar{f}) \qquad}&
\end{align}
\end{coro}
\begin{proof}
Let $f\in\hom_{\Ab^\Delta}(QA,QB)$ and $j,\pe$ and $h$ be as in the theorem.
Define $\bar{f}:=\pe N(f)j$. Because $\pe j=1$, $\overline{Qg}=g$. Using the naturality of $j$ and $\pe$, one checks further that $f\mapsto\bar{f}$ is natural.
Let $\delta=N(f)-NQ(\bar{f})$ and put
$$
\kappa=\kappa_f:=h\delta+\delta h-[h\delta,\partial]h
$$
One checks that $[\kappa,\partial]=\delta$, whence
$H_f:=(Nf,\kappa,NQ\bar{f})$ is a homomorphism $NQA\to PNQB$ with
$\epsilon H_f=(Nf,NQ\bar{f})$. The naturality of $H:f\mapsto H_f$ follows from
that of $h$.
\end{proof}
\begin{powers}\label{papa}
Let $A\in\Ab^\Delta$, $X\in\sets^{\Delta^{op}}$. Put
\begin{equation}
(A^X)^n:=\prod_{x\in X_n}A^n
\end{equation}
If $\alpha\in\hom_{\Delta}([n],[m])$ and $a\in (A^X)^n$, define
$\alpha(a)_x=\alpha(a_{\alpha x})$ $(x\in X_m)$.
The dual $\Z[X]^{\vee}:[n]\mapsto \hom_\Z(\Z[X_n],\Z)$ of
the simplicial free abelian group $\Z[X]$ is a cosimplicial
group. Consider the cosimplicial tensor product
$A\otimes\Z[X]^{\vee}:[n]\mapsto A^n\otimes\Z[X_n]^{\vee}$. There is a natural
homomorphism
\begin{equation}
\eta:A\otimes\Z[X]^{\vee}\to A^X,\qquad \eta(a\otimes \phi)_x=a\phi(x)
\end{equation}
In case each $X_n$ is finite, $\eta$ is an isomorphism.
Dualizing the statement in \cite{mac} --next after 8.9-- we get that
the composite of the normalized shuffle map
$NA\otimes N\Z[X]^{\vee}\to N(A\otimes \Z[X]^{\vee})$  with Alexander-Whitney map
$N(A\otimes \Z[X]^{\vee})\to NA\otimes N\Z[X]^{\vee}$ is the identity. Thus
$NA\otimes N\Z[X]^{\vee}$ is a deformation retract of $N(A\otimes \Z[X]^{\vee})$.
In particular $PNA=NA\otimes N\Z[\Delta[1]]^{\vee}$ is a
deformation retract of $N(A^{\Delta[1]})$. Recall two cosimplicial maps
$f_0,f_1:A\to B$ are called
{\it homotopic} if $(f_0,f_1):A\to B\times B=B^{\Delta[0]\coprod\Delta[0]}$
can be lifted to a map $H:A\to B^{\Delta[1]}$. From what we have just
seen it is clear that $f_0,f_1$ are homotopic in this sense if and only
if $Nf_0$,$Nf_1$ are cochain homotopic. (The dual of this assertion is proved
in \cite{dold}.)
Let $\C$ be either of $\cha$, $\Ab^\Delta$. We write $[\C]$
for the category with the same objects as $\C$, but where the
homomorphisms are the homotopy classes of maps in $\C$.
\end{powers}
\begin{propo}
The functor $Q$ induces an equivalence of categories $[\cha]\to[\Ab^\Delta]$.
\end{propo}
\begin{proof}
If $A\in \Ab^\Delta$, then $A=KNA$. By Theorem \ref{kequivq}, $NA$ is homotopy
equivalent to $NQA$. Thus $A$ is homotopy equivalent to $KNQA=QA$. It remains
to show that the following map is a bijection
\begin{equation*}
[Q]:\hom_{[\cha]}(A,B)\to\hom_{[\Ab^\Delta]}(QA,QB).
\end{equation*}
It is clear from the previous corollary that
the composite of $Q$ with
\begin{equation}\label{nhom}
[N]:\hom_{[\Ab^\Delta]}(QA,QB)\to \hom_{[\cha]}(NQA,NQB)
\end{equation}
is a bijection. But \eqref{nhom} is
bijective by \ref{papa}.
\end{proof}

\begin{defi}
Give $\cha$ the closed model category structure in which a map is a {\it fibration} if
it is surjective codimensionwise, a {\it weak equivalence} if it is a
quasi-isomorphism, and a cofibration if it has Quillen's left lifting
property ($LLP$, see \cite{qui}) with respect to those fibrations which
are also weak equivalences ({\it trivial fibrations}). All this structure
carries over to $\Ab^\Delta$ using the category equivalence
$N:\Ab^\Delta\to\cha$.
In the lemma below $RLP$ stands for right lifting property in the sense
of \cite{qui}.
\end{defi}

\bigskip

\begin{nota}
In the next lemma and further below, we use the following notation.
If $n\ge 0$, we write $\Z<n,n+1>$ for the mapping cone of the identity map
$\Z[n]\to \Z[n]$.
\end{nota}
\begin{lemm}\label{equifib}
Let $f:E\to B$ be a homomorphism of cosimplicial abelian groups.
We have:
\smallskip

\noindent{i)} $f$ is a fibration if and only if for all $n\ge 1$ $f$ has the
$RLP$ with respect to $0\to Q\Z<n-1,n>$.

\smallskip

\noindent{ii)} $f$ is a trivial fibration if and only if for all $n\ge 1$ $f$
has the $RLP$ with respect to the natural inclusion
$Q\Z[n]\hookrightarrow Q\Z<n-1,n>$.

\end{lemm}
\begin{proof}
Let $f:C\to D$ be a cochain map. By the theorem, $Kf$ is a retract of $Qf$.
Thus every map having the $RLP$ with respect to $Qf$ also has it with respect
to $Kf$. The lemma follows from this applied to the cochain maps
$0\to \Z<n-1,n>$ and $\Z[n]\hookrightarrow \Z<n-1,n>$.
\end{proof}


\section{Monoidal structure}\label{Qismonoidal} 

Consider the map $\theta:TV^*\to TV^*$,
\begin{equation}\label{teta}
\theta(v_i)=v_i^2,\qquad \theta(xy)=\theta(x)y+(-1)^{|x|}x\theta(y).
\end{equation}
The second identity says that $\theta$ is a homogeneous derivation of degree $+1$. 
Note it follows from \eqref{teta} that $\theta^2 = 0$.

\begin{lemm}
\label{alphatheta}
For every $\alpha \in Map([n], [m])$ and $x \in TV^n$,  $[\alpha, \theta](x) = [v_{\alpha(0)}, \alpha(x) ]$.
\end{lemm}
\begin{proof}
Both sides of the identity we have to prove are derivations. Thus it suffices to show they agree on the generators $v_i$, and this is straightforward.
\end{proof}

\begin{them}
\label{Qisstrong}
Let $A, B \in \cha$ and $\theta$ as defined in \eqref{teta} above. Consider the tensor product of $\fin$-abelian groups $QA \otimes QB : [n] \to Q^nA \otimes Q^nB$. 
The map $\upsilon:QA \otimes QB \to Q(A \otimes B)$ given 
by the following formula is an isomorphism in $\Ab^\fin$, and makes 
$Q : \cha \to \Ab^\fin $ a strong monoidal functor: 

\[
\upsilon((a \otimes x) \otimes (b \otimes y)) = a \otimes b \otimes xy + (-1)^{\vert a \vert}a \otimes db \otimes \theta(x)y
\]
\end{them}
\begin{proof}
It is clear that the following map is an isomorphism of abelian groups: 
\[
g :(a \otimes x) \otimes (b \otimes y) \mapsto (a\otimes x)\cdot(b\otimes y):
=a \otimes b \otimes xy
\]
Because $h := \upsilon- g$ is homogeneous of 
degree +1 and $h^2 = 0$, $\upsilon$ is a group isomorphism. 
That $\upsilon$ is a homomorphism in $\Ab^{\fin}$ follows straightforwardly using 
Lemma \ref{alphatheta}. In order to see that $Q$ is strong monoidal, we must check 
that the two diagrams involving the unit object of $\Ab^\fin$ commute,
which is immediate, and also the following associativity condition
for $\alpha\in QA$, $\beta\in QB$ and $\gamma\in QC$
\begin{equation}\label{aso1}
\upsilon(\upsilon(\alpha\otimes\beta)\otimes\gamma)=\upsilon(\alpha\otimes \upsilon(\beta\otimes\gamma))
\end{equation}
Writing this in terms of $g$ and $h$, and because $g$ is associative, we obtain
\begin{align}\label{aso2}
h(h(\alpha\otimes\beta)\otimes\gamma)-h(\alpha\otimes h(\beta\otimes\gamma))=
&\alpha\cdot h(\beta\otimes\gamma))-h(\alpha\cdot\beta\otimes\gamma)+\\
h(\alpha\otimes \beta\cdot\gamma)-h(\alpha\otimes\beta)\cdot\gamma\nonumber
\end{align}
For $\alpha=a\otimes x$, $\beta=b\otimes y$ and $\gamma=c\otimes z$,
the left hand side of \eqref{aso2} is
\[
(-1)^{|y|+1}adbdc\otimes\theta(\theta(x)y)z+(-1)^{|x|+|y|+1}adbdc
\otimes\theta(x)\theta(y)z=0
\]
This is zero because $\theta$ is a square-zero derivation. Thus \eqref{aso2}
says that $h$ is a Hochschild $2$-cocycle, which follows from the fact
that both $d$ and $\theta$ are derivations.
\end{proof}

\begin{product}
Let $A \in DGR^*$, $m : A \otimes A \to A$ the multiplication map. Consider the composite 
\[
\circ : 
\xymatrix{
QA \otimes QA \ar[r]^{\upsilon} & Q(A \otimes A) \ar[r]^{\quad Qm}& QA.
} 
\] 
We have
\begin{equation}\label{fedo}
(\omega\otimes x)\circ (\eta\otimes y):=\omega\eta\otimes xy+(-1)^{|x|}
\omega d\eta\otimes\theta(x)y
\end{equation}
By construction, $(QA, \circ)$ is a $\fin$-ring. Note that each term $\F_nQA$ of the
filtration \eqref{filq} is a $\fin$-ideal. The associated graded $\fin$-ring
is $A\boxtimes TV$ equipped with the product inherited from
$A\boxtimes TV\subset A\otimes TV$.
Thus we may view $QA$ as a deformation of $A\boxtimes TV$.
One checks that the kernel of the map $\pe:QA\fib KA$ of \eqref{mape} is
an ideal for $\circ$. Hence $KA$ inherits a $\fin$-ring structure; using
the definition of $\theta$ we get that the induced product on
$KA=A\boxtimes\Lambda V$ is just that coming from $A\otimes \Lambda V$:
\begin{equation}\label{trinus}
(a\otimes x)(b\otimes y)=ab\otimes x\land y.
\end{equation}
\end{product}

\section{Comparison with the shuffle product}\label{comshu}
Let $R$ be a simplicial ring. Consider the direct sum of its homotopy groups
\begin{equation}\label{defpi}
\pi R:=\bigoplus_{n=0}^\infty\pi_n R.
\end{equation}
Recall that the shuffle product $\star$ makes $\pi R$ into a graded ring. If
moreover $R$ is a $\fin$-ring, then the Connes operator
$B:\pi_* R\to \pi_{*+1} R$ is a derivation, so that
$\pi R=(\pi R,\star,B)$ becomes
in fact a $DG$-ring. This follows from the version of \cite{lod}, 4.3.3. for cyclic modules, the same which is used without further proof in \cite{lod}, 4.3.7-8. Hence we have a functor
\begin{equation}\label{funpi}
\ass^\fin\to DGR^*,\qquad R\mapsto \pi R.
\end{equation}

\begin{propo}\label{propi}
Let $A\in DGR^*$. Consider the natural isomorphism of graded abelian groups
induced by the map $l$ of {\rm\ref{kequivq} ii)} 
\begin{equation}\label{liso}
l:A\weq \pi QA.
\end{equation}
The map \eqref{liso} is an isomorphism of $DG$-rings.
In particular the functor \eqref{funpi} is a left inverse of $Q$.
\end{propo}
\begin{proof}
By \ref{kequivq}, $l$ induces a cochain isomorphism $(A,d)\cong (\pi QA,B)$.
It remains to show that the induced map is a ring homomorphism.
Recall the formula for the shuffle product $\star$ involves
degeneracies and shuffles. Keeping in mind that the degeneracies in $QA$ are
of the form $s_i=1\otimes\partial_{i+1}$ with $\partial_j$ the coface of
$TV$, we get the following identity for $a\in A^n$, $b\in A^m$:
\begin{align*}
l(a)\star l(b)=&(a\otimes \epsilon_n)\star(b\otimes\epsilon_m)\\
\equiv&ab\otimes\epsilon_{n}\star\epsilon_{m}\mod N\F_{n+m+1}Q^{n+m}A\\
=&ab\otimes\epsilon_{n+m}=l(ab)
\end{align*}
This finishes the proof, since $\pi_{n+m} N\F_{n+m+1}QA=0$ by the proof
of \ref{kequivq}.
\end{proof}

\section{Noncommutative Hochschild-Kostant-Rosenberg\\
         and Loday-Quillen theorems}\label{nchkr}

Recall from \cite{we} that for every algebra $S$ over a commutative ring $R$
which is central in $S$
there is defined a cyclic $R$-module $C_*(S/R)$.
Recall also that the normalization of $C_*(S/R)$ is the mixed complex of 
noncommutative
differential forms \cite{cq2} $NC_*(S/R)=\Omega_RS$.
The Hochschild-Kostant-Rosenberg theorem (\cite{we}, Ex. 9.4.2)
says that if $R$ and $S$ are commutative, $R$ noetherian, and $R\to S$ an
essentially of finite type,
smooth homomorphism, then the canonical
map from commutative differential forms to Hochschild homology
induced by the shuffle product is an isomorphism:
\begin{equation}\label{hkr}
\Omega_{S/R}^*=\Lambda^* HH_1(S/R)\weq HH_*(S/R).
\end{equation}
If $R\supset\Q$ the inverse of \eqref{hkr} is induced by the
homomorphism
\begin{equation}\label{pecom}
\Omega_RS\to \Omega_{S/R} \qquad a_0da_1\dots da_n\mapsto
\frac{1}{n!}a_0da_1\land\dots\land da_n
\end{equation}
Here the boundary operators are the Hochschild boundary $b$ on $\Omega_RS$
and the trivial boundary on $\Omega_{S/R}$.
Moreover, as \eqref{pecom} maps $B$ to $d$, it is in fact a mixed
complex equivalence
\begin{equation*}
(\Omega_RS,b,B) \weq (\Omega_{S/R},0,d)
\end{equation*}
We will prove an analogue of this which holds for not necessarily commutative
$R$ and $S$.
Note that if $R$ and $S$ are commutative then $C_*(S/R)$ is just the coproduct $\fin$-algebra $\bigotimes_RS$ considered as a cyclic module. 
The analogue concerns
the coproduct $\fin$-ring $\coprod_RS$ which arises from a ring homomorphism
$R\to S$ of not necessarily commutative rings. We show in \ref{lq} below
that
there is an equivalence of mixed complexes $(\Omega_RS,0,d)\to(N\coprod_RS,\mu,B)$, valid without restrictions on the characteristic. We deduce
this from \ref{kequivq} and from \ref{qomega} below, where
we show that $\coprod_RS=Q\Omega_RS$. In particular the isomorphism
$\Omega^*_RS\cong HH_*(N\coprod_RS,\mu,B)=\pi_*\coprod_RS$ is \eqref{liso},
which is a ring homomorphism for the product of forms and the
shuffle product (by \ref{propi}) just like the Hochschild-Kostant-Rosenberg
isomorphism \eqref{hkr}. Note further the analogy between \eqref{pecom}
and the rescaled map $\pe$ of \ref{rescape}.

To prove the isomorphism $Q\Omega_RS\cong\coprod_RS$ we show first that
$Q$ has a right adjoint (\ref{adjoint}). In the next lemma we use the
symbol $T$ for both the the tensor $\fin$- and $DG$- rings.

\begin{lemm}\label{qttq}
Let $(U,d)\in \cha$. Then there is a natural isomorphism of $\fin$-rings
$TQU\iso QTU$.
\end{lemm}
\begin{proof}
This is a formal consequence of Theorem \ref{Qisstrong}.
\end{proof}

\begin{nota} The following $DG$-rings shall be considered often in what
follows
\begin{equation}
D(n):=T\Z<n,n+1>\supset S(n):=T\Z[n]
\end{equation}
\end{nota}

\begin{coro}\label{corocoprodn}
Let $I$ be a set, $n_i\ge 0$. Then
$$Q(\coprod_{i\in I}D(n_i))=\coprod_{i\in I}QD(n_i).$$
\end{coro}
\begin{proof}
\begin{align*}
Q(\coprod_{i\in I}D(n_i))=&Q(\coprod_{i\in I}T(\Z<n_i,n_{i+1}>))=QT(\bigoplus_{i\in I}\Z<n_i,n_{i+1}>)\\
=&TQ(\bigoplus_{i\in I}\Z<n_i,n_{i+1}>)=T(\bigoplus_{i\in I}Q\Z<n_i,n_{i+1}>)\\
=&\coprod_{i\in I}TQ\Z<n_i,n_{i+1}>=\coprod_{i\in I}QD(n_i).
\end{align*}
\end{proof}
\begin{propo}\label{adjoint}
Let $\ass$ be the category of associative unital rings and  $DGR^*$ that
of cochain differential graded rings. The functor $Q:DGR^*\to \ass^\fin$
has a right adjoint.
\end{propo}
\begin{proof}
This is an adaptation of the proof of the dual of Freyd's Special Adjoint
Theorem (\cite{cwm}, Chapter V,\S8, Theorem 2).
Let $B\in\ass^\fin$. Put
\begin{equation}
DB:=\coprod_{n\ge 0}\coprod_{\hom(QD(n),B)}D(n)
\end{equation}
If $s\in \hom(QD(n),B)$, write $j_s:D(n)\to DB$ for the corresponding
inclusion. Define $\alpha:QDB\to B$ by $\alpha j_s=s$. Consider the
two-sided $\fin$-ideal
\begin{equation}
DB\triangleright K:=\sum\{I\triangleleft DB:\alpha(QI)=0\}
\end{equation}
Set $PB:=DB/K$. Because $Q:\cha\to \Ab^\fin$ is exact, we have a natural map
$\hat{\alpha}$ making the following diagram commute
\begin{equation}\label{alphabar}
\xymatrix{QDB\ar[r]^\alpha\ar[d]&B\\
          QPB\ar[ur]_{\hat{\alpha}}&}
\end{equation}
Hence $(PB,\hat{\alpha})$ is an object of the category $Q\uparrow B$
(notation is as in \cite{cwm}). We
shall see it is final, which proves that $P$ is right adjoint to $Q$.
Let $(R,f)\in Q\uparrow B$. Put
\begin{equation*}
ER:=\coprod_{n\ge 0}\coprod_{\hom(D(n),R)}D(n).
\end{equation*}
If $r:D(n)\to R$ is a homomorphism, write $i_r:D(n)\to ER$
for the corresponding inclusion. Consider the homomorphisms
$\pi:ER\to R$, $\pi i_r=r$ and $g:ER\to DB$, $g i_r=j_{fQr}$.
We claim that the following diagram commutes
\begin{equation}\label{disdiag}
\xymatrix{QER\ar[r]^{Qg}\ar[d]_{Q\pi}&QDB\ar[d]^\alpha\\
            QR\ar[r]_f&B}
\end{equation}
Indeed by \ref{corocoprodn}, commutativity can be checked at each ``cell''
$Q(D(n))$ where it is clear.
Using \eqref{disdiag} together with the exactness of $Q$, we get that
$g(\ker\pi)\subset K$. Thus $g$ induces a map $\hat{g}$ making the
following diagram commute
\begin{equation}
\xymatrix{ER\ar[r]^g\ar[d]_\pi&DB\ar[d]\\
           R\ar[r]_{\hat{g}}&PB}
\end{equation}
It follows that also the following commutes
\begin{equation}\label{alsoco}
\xymatrix{QER\ar[r]^{Qg}\ar[d]_{Q\pi}&QDB\ar[dr]^\alpha&\\
           QR\ar[r]_{Q\hat{g}}&QPB\ar[r]_{\hat{\alpha}}&B}
\end{equation}
Putting together the latter diagram with \eqref{alphabar} and \eqref{disdiag}
we get that $fQ(\pi)=\hat{\alpha}Q(\hat{g})Q(\pi)$. Because $\pi$ is surjective
and $Q$ exact, we conclude $f=\hat{\alpha}Q(\hat{g})$; in other words
$\hat{g}$ is a homomorphism
$(R,f)\to (PB,\hat{\alpha})$ in $Q\uparrow B$. Let $h:(R,f)\to (PB,\hat{\alpha})$ be another. Lift $\hat{h}$ to a map $h:ER\to DB$. Then by \eqref{alsoco},
\begin{equation*}
\alpha Q(h)=\hat{\alpha}Q(\hat{h}\pi)=fQ(\pi)=\alpha Q(g)
\end{equation*}
Hence the image of $g-h$ lands in $K$, and therefore $\hat{g}=\hat{h}$.
\end{proof}

\begin{rema}
Essentially the same proof as that of the Theorem above shows that
also $Q:\cha\to\Ab^\fin$ has a right adjoint. One just has to replace
$\coprod$ and $D(n)$ for $\oplus$ and $\Z<n,n+1>$.
\end{rema}

\begin{them}\label{qomega}
Let $R\to S$ be a ring homomorphism, $R\uparrow\ass$ the category
of $R$-algebras, $\coprod_R$ the coproduct in $R\uparrow\ass$,
$\coprod_RS$ the $\fin$-ring of Section {\rm \ref{linear}} above and
$\Omega_RS$ the $R$-$DG$-algebra of relative noncommutative
differential forms of \cite{cq}. Then
$Q(\Omega_RS)=\coprod_RS$.
\end{them}
\begin{proof}
The $\fin$-ring $\coprod_RS$ is characterized by the following property
\begin{equation}
\hom_{(R\uparrow\ass)^\fin}({\coprod}_{R}S,C)=\hom_{R\uparrow\ass}(S,C^0)
\end{equation}
We must show $Q\Omega_RS$ has the same property. On the other hand
we have
\begin{equation}\label{omegaab}
\hom_{R\uparrow DGR^*}(\Omega_RS,X)=\hom_{R\uparrow\ass}(S,X^0)
\end{equation}
Here we identify $R$ with the $DGR^*$ concentrated in codimension $0$ with
trivial derivation. Let $P$ be the right
adjoint of $Q:DGR^*\to\ass^{\fin}$; its existence is guaranteed by Proposition
\ref{adjoint}. Identifying $R$ with the constant
$\fin$-ring, noting that $QR=R$ and using \eqref{omegaab},
we obtain
\begin{align*}
\hom_{R\uparrow(\ass^\fin)}(Q(\Omega_{R}S),C)&=\hom_{R\uparrow DGR^*}(\Omega_RS,PC)\\
&=\hom_{R\uparrow \ass}(S,PC^0)
\end{align*}
Therefore to prove the corollary it suffices to show that $PC^0=C^0$.
We have
\begin{align}\label{qz01pc}
PC^0=&\hom_{\cha}(\Z<0,1>,PC)\\
=&\hom_{DGR^*}(D(0),PC)\nonumber\\
                            =&\hom_{\ass^\fin}(QD(0),C)\nonumber\\
                            =&\hom_{\ass^\fin}(TQ\Z<0,1>,C)\text{ (by \ref{qttq})}\nonumber\\
                            =&\hom_{\Ab^\fin}(Q\Z<0,1>,C)\nonumber\\
\end{align}
By definition
\begin{equation}
Q^n\Z<0,1>=\Z<0,1>\boxtimes TV^n=
\Z(1\otimes 1)\oplus\bigoplus_{i=1}^n\Z(1\otimes v_i)
\end{equation}
Put $e_0=1\otimes 1$, $e_i=1\otimes v_i+e_0$ $1\le i\le n$. It follows
from \eqref{alphav} that
 $\alpha(e_i)=e_{\alpha (i)}$ for all $\alpha:[n]\to [m]\in\fin$. Therefore
$Q\Z<0,1>=\bigoplus \Z$, whence \eqref{qz01pc} equals
\begin{align*}
=&\hom_{\Ab^\fin}(\bigoplus \Z,C)=C^0 \qquad \qquad
\end{align*}
\end{proof}

\begin{coro}
\label{lq}
View the $\fin$-ring $\coprod_RS$ as a cyclic module by
restriction, and consider its associated normalized mixed complex
$(N\coprod_RS,\mu,B)$. Then the map $l$ of Theorem {\rm\ref{kequivq}} is
a mixed complex equivalence
$l:(\Omega_RS,0,d)\to(N\coprod_RS,\mu,B)$
\end{coro}

\bigskip

\begin{rema}
As a particular case of Theorem \ref{qomega} we get a ring isomorphism
\begin{equation}
S\coprod_RS\cong Q^1\Omega_RS=\Omega_RS\boxtimes TV^1\cong\Omega_RS
\end{equation}
Here $\Omega_RS$ is equipped with the product $\circ$ of \eqref{fedo}. A
similar isomorphism but with a different choice of $\circ$ was proved
by Cuntz and Quillen in \cite{cq} Proposition 1.3, under the stated assumption
that $R=\mathbb{C}$. Their choice of $\circ$
actually works whenever $2$ is invertible in $R$, and
the rings which arise from $\Omega_RS$ with our product and that of \cite{cq}
are isomorphic in that case. Hence \ref{qomega} may be viewed as a strong
generalization
of Cuntz-Quillen's result.
\end{rema}

\section{Comparison with Nuss' product}\label{conu}

In \cite{nuss}, P. Nuss considers the
``twist''
$$
\tau:S\otimes_RS\to S\otimes_RS, \qquad \tau(s\otimes t)=st\otimes 1+1\otimes
st-s\otimes t
$$
It is clear that $\tau^2=1$ and that, for the multiplication map $\mu_0:S\otimes_RS\to S$, we have $\mu_0\tau=\mu_0$. He shows further (\cite{nuss}, 1.3) that
$\tau$ satisfies the Yang-Baxter equation. Using $\tau$, he introduces a
ring structure on the $n+1$ fold tensor power $S\otimes_R\dots\otimes_RS$
for all $n\ge 1$, by a standard procedure (use Proposition 2.3 of \cite{who} and
induction).
We want to reinterpret this product in a different way. For
this consider the (Amitsur) cosimplicial $R$-bimodule
\begin{equation*}
{\bigotimes}_R S
:[n]\mapsto {\bigotimes_{i=0}^n}{}_{{}_R} S
\end{equation*}
By definition of $\Omega_RS$, we have $N({\bigotimes}_R S)=\Omega_RS$. Hence
the Dold-Kan correspondence gives an isomorphism of cosimplicial $R$-bimodules
\begin{equation}\label{isonus}
{\bigotimes}_R S\cong K\Omega_RS
\end{equation}
On the right hand side we also have the product \eqref{trinus}. It is noted
in \cite{nuss} that \eqref{isonus} is a ring isomorphism in codimension
$\le 1$.
The next Proposition shows it is actually a ring isomorphism in all 
codimensions.

\begin{propo}\label{nusbo}
Equip ${\bigotimes_R}S$ with the product defined in \cite{nuss} and
$K\Omega_RS$ with that given by \eqref{trinus}. Then {\rm\eqref{isonus}} is
an isomorphism of $\fin$-rings.
\end{propo}
\begin{proof}
Write $\bullet$ for Nuss' product.
Consider the following map
\begin{equation*}
\delta_i=\partial_{i+1}^{n-i}\partial_0^i\in\hom_\Delta([0],[n]) \qquad
(0\le i\le n).
\end{equation*}
One checks the following identities hold in $\bigotimes^n_RS$, for $a, b \in S$:
\begin{equation}\label{formudelta}
\delta_i(a)\bullet\delta_j(b)=\begin{cases}
\partial_{j+1}^{n-j}\partial_{i+1}^{j-i-1}\partial_0^i(a\otimes b)& i<j\\
\delta_i(ab) & i=j \\
-\delta_j(a)\bullet\delta_i(b)+\delta_i(ab)+\delta_j(ab) &i>j
\end{cases}
\end{equation}
In particular $\delta_i:S\to\bigotimes^n_RS$ is a ring
homomorphism for $\bullet$. By the universal property of
$\coprod_R^nS$, we have a unique ring homomorphism
$\alpha^n:\coprod_R^nS\to \bigotimes^n_RS$ satisfying
$\alpha^n\delta_i=\delta_i$ for all $i$. By \eqref{formudelta},
\begin{align*}
s_0\otimes\dots\otimes s_n=&\delta_0(s_0)\bullet\dots\bullet\delta_n(s_n)\\
                          =&\alpha(\delta_0(s_0)\dots\delta_n(s_n)).
\end{align*}
Thus $\alpha$ is surjective. On the other hand the composite of $\alpha$
with the isomorphism
$Q\Omega_RS\weq\coprod_RS$ sends $ds\otimes v_i$ to
$q_i(s):=\delta_i(s)-\delta_0(s)$. But it follows from \eqref{formudelta}
that
\begin{equation}
q_i(a) \bullet q_j(b)=\begin{cases}-q_j(a) \bullet q_i(b)& (i\ne j)\\
                           0 & (i=j)\end{cases}
\end{equation}
Thus $\alpha$ descends to a ring homomorphism
$\overline{\alpha}:K\Omega_RS\to{\bigotimes}_RS$. On the other hand we
have an $R$-linear map $\beta:\bigotimes_RS\to K\Omega_RS$,
$\beta(s_0\otimes\dots\otimes s_n)=
\delta_0(s_0)\dots\delta_n(s_n)$. Clearly $\alpha\beta=1$. To finish the proof
it suffices to show that $\beta$ is surjective. But we have
\begin{align*}
&a_0da_1\dots da_r\otimes v_{i_1} \wedge \dots \wedge v_{i_r}\\
= &\delta_0(a_0)(\delta_1(a_1)-\delta_0(a_1))\dots(\delta_r(a_r)-\delta_0(a_r))\\
\equiv &\delta_0(a_0)\dots\delta_r(a_r)\mod\bigoplus_{i=0}^{r-1}\Omega^i_RS\otimes\Lambda^iV\\
= &\beta(a_0\otimes\dots\otimes a_r\otimes 1 \otimes \dots \otimes 1).
\end{align*}
Hence it follows by induction on $r$, that $\Omega_R^rS \otimes \Lambda^rV$ is included in the image of $\beta$.
\end{proof}

\section{Dold-Kan equivalence for rings}
\label{equimain}
\begin{defi}\label{defcm}
Let $f:R\to S$ be a homomorphism in $DGR^*$. We say that $f$ is a
{\it weak equivalence} if it induces an isomorphism in cohomology.
We call $f$ a {\it fibration} if each $f^n:R^n\to S^n$ is surjective, and a {\it cofibration}
if it has the left lifting property ($LLP$) of \cite{qui} 
 with respect to those fibrations which are also weak
equivalences (trivial fibrations). Similarly, a map $g:A\to B$ of cosimplicial
rings is a {\it weak equivalence} if it induces an isomorphism in cohomotopy,
a {\it fibration} if each $g^n:A^n\fib B^n$ is surjective and a {\it cofibration} if it has the $LLP$ with respect to trivial fibrations. 
It is proved
in \cite{jar} that the structure just defined makes $DGR^*$ closed model. 
The next proposition shows that the same is valid for cosimplicial rings.
\end{defi}
\begin{propo}\label{yescm}
With the notions of fibration, cofibration and weak equivalence defined
in {\rm\ref{defcm},} $\ass^\Delta$ is a closed model category.
\end{propo}
\begin{proof}
A commutative version of this is given in \cite{to}, Theorem 2.1.2.
Essentially the same proof works in the noncommutative case; simply
substitute the coproduct $\coprod$ of $\ass$ for $\otimes$, which
is the coproduct in the category $\comm$ of commutative rings. One only
has to check that for all $n\ge 0$, the structure maps 
$\Z\to D(n):=TK\Z<n,n+1>\in\ass^\Delta$ induce weak equivalences
\begin{align}
A\weq &A \coprod D(n)\qquad (A\in \ass^\Delta)\label{cscm}
\end{align}
For this we imitate Jardine's argument (\cite{jar}).
We observe that if $A \in \ass^\Delta$ and we write $C(n) = K\Z<n,n+1>$ then 
there is an isomorphism of cosimplicial groups
\begin{align*}
A \coprod D(n) = & A[C(n)] := \\
&A \oplus (A \otimes C(n) \otimes A) \oplus (A \otimes  C(n) \otimes A \otimes  C(n) \otimes A) \oplus ...
\end{align*}
with the product defined by 
\begin{align*}
(a_1 \otimes c_1 \otimes a_2 \otimes ... \otimes c_k \otimes a_{k+1})(a'_1 \otimes c'_1 \otimes a'_2 \otimes ... \otimes c'_l \otimes a'_{l+1}) = & \\
(a_1 \otimes c_1 \otimes ... \otimes c_k  \otimes a_{k+1}a'_1 \otimes ... \otimes c'_l \otimes a'_{l+1})
\end{align*}
and cofaces and codegeneracies induced by those of $A$ and $C(n)$.
Thus to prove \eqref{cscm} it suffices to show that if $C$ and $D$ are
cosimplicial groups and $D$ is contractible, then the inclusion $\iota:C \to C[D]$ is 
a quasi-isomorphism. But $\coker\iota$ is a sum of cosimplicial groups each
of which is isomorphic to one of the form $C \otimes D \otimes \dots \otimes D \otimes C$.
Hence it suffices to show that $D\otimes D'$ is contractible if $D$ is. This
latter statement follows from the following property of the cosimplicial
path functor (see \cite{to}, page 30):
$$D^{\Delta[1]}\otimes D'= (D\otimes D')^{\Delta[1]}.$$
\end{proof}

\begin{lemm}\label{qpreserve}
\item{i)}The functor $\tilde{Q}:DGR^*\overset{Q}\longrightarrow \ass^\fin\overset{forget}\longrightarrow\ass^\Delta$ preserves colimits, finite limits, 
cofibrations, fibrations, and weak equivalences.
\item{ii)} Let $K:DGR^*\to \ass^\Delta$ be the functor sending
$A\mapsto KA$ where $KA$ is equipped with the product \eqref{trinus}.
Then there is a natural isomorphism of left derived functors 
$\mathbb{L}\tilde{Q}\iso\mathbb{L}K$.
\end{lemm}
\begin{proof}
Limits and colimits in $\ass^\Delta$ are computed codimensionwise, and the
same is true in $\ass^\fin$. In particular the forgetful functor preserves
limits and colimits. The functor $Q:DGR^*\to\ass^\fin$ preserves colimits
by Proposition \ref{adjoint}. Thus $\tilde{Q}$ preserves colimits. On the other
hand limits in $\ass^\fin$ can be computed in $\Ab^\fin$. As 
$Q:\cha\to\Ab^\fin$ is exact and preserves direct sums, it follows that 
$\tilde{Q}$ preserves finite limits. Similarly, as the forgetful functors
$DGR^*\to \cha$ and $\ass^\Delta\to\Ab^\Delta$ as well as $Q:\cha\to \Ab^\fin$ 
preserve fibrations and weak equivalences, it follows that $\tilde{Q}$ does.
One checks, using Lemma \ref{qttq}, that $\tilde{Q}$ preserves the basic cofibrations $S(m)\to D(m)$, $\Z\to D(m)$. Because it also preserves colimits it
follows that if $m_i,i\in I$ is a family of positive integers and 
$e_i:S(m_i)\to X$ $(i\in I)$ a family of maps, then the
following maps are cofibrations:

\begin{align*}
\tilde{Q}(X\cofi &X{\coprod}_{\coprod_{i\in I}S(m_i)}\coprod_{i\in I} D(m_i))\\
\tilde{Q}(X\cofi &X\coprod\coprod_{i\in I} D(m_i))
\end{align*}
But by the proof of \ref{yescm} and  the remark on page 23 of \cite{bogu},
every cofibration in $DGR^*$ is a retract of one obtained as a colimit
of such cofibrations. Hence $\tilde{Q}$ preserves all cofibrations.
Thus i) is proved. As shown in Section \ref{Qismonoidal}, the natural weak equivalence 
$\pe:\tilde{Q}A \to KA$ of \ref{kequivq} is a homomorphism
of cosimplicial rings. This proves ii).
\end{proof}

\bigskip
\begin{rema}\label{LneQ}
A functor $L_*$ with properties similar to those proved for $\tilde{Q}$ in Lemma
\ref{qpreserve} is considered in \cite{ship} for the dual situation of
chain $DG$- and simplicial rings. The authors use the shuffle product to
make the normalized chain complex of a simplicial ring into a chain $DG$-ring,
thus obtaining a functor $N_*:\ass^{\Delta^{op}}\to DGR_*$. The functor
$L_*$ is defined as the left adjoint of $N_*$. Dually, one can equip
the normalized complex of a cosimplicial ring with the shuffle product,
consider the resulting functor $\tilde{N}:\ass^\Delta\to DGR^*$ and take
its left adjoint $L^*$. However we point out that $L^*$ and $\tilde{Q}$ are not
isomorphic. In other words $\tilde{Q}$ is not left adjoint to $\tilde{N}$.
To see this, note that, by \ref{qttq}, if $A\in \cha$, then  
$\hom_{\ass^\Delta}(\tilde{Q}TA, R)=\hom_{\Ab^\Delta}(QA,R)$, while
$\hom_{DGR^*}(TA,\tilde{N}R)=\hom_{\cha}(A,NR)=\hom_{\Ab^\Delta}(KA,R)$. 
Hence if $\tilde{Q}$ were left adjoint to $\tilde{N}$, then $K$ and $Q$ should
be isomorphic as functors $\cha\to\Ab^\Delta$, 
which is clearly false.
\end{rema}

\bigskip

\begin{rema}\label{bp}
We have seen in Proposition \ref{adjoint} that $Q$ has a right adjoint $P$. 
Since the forgetful functor $U : \ass^\fin \to \ass^\Delta$ also has a right 
adjoint (\cite{cwm}, X.3.2), and $\bq=UQ$, it follows that $\bq$ is the left 
adjoint of an adjoint pair $(\bq,\bp)$. On the other hand, by lemma 
\ref{qpreserve}.i), we know that $\bq$ preserves cofibrations and weak 
equivalences, and thus it is the left adjoint of a Quillen adjoint functor 
pair (\cite{hov}, def. 1.3.1).
\end{rema}

\begin{them}
\label{main1}
The adjoint functors $\bq : DGR^* \leftrightarrows \ass^\Delta : \bp$
of {\rm \ref{qpreserve} i)} and {\rm\ref{bp}} form a Quillen equivalence
in the sense of {\rm \cite{hov} 1.3.12.}\end{them}
\begin{proof}
Let $g: R := \bp(S)^c \trfib \bp(S)$ be the functorial cofibrant replacement
obtained by the small object argument.
Since the functor $\bq$ reflects weak equivalences, it suffices to show that 
the adjoint map $f: \bq R \to S$ is a weak equivalence 
(\cite{hov}, Theorem 1.3.16.). 
We note for future use that by the small object argument and because $\bq$ and $\bp$ are adjoint, the dotted arrow in the diagram below exists whenever the top horizontal arrow
is in the image of $\tilde{Q}:\hom_{DGR^*}(S(m),R)\to \hom_{Rings^\Delta}(\tilde{Q}S(m),\tilde{Q}R)$.
\begin{equation}\label{liftexi}
\xymatrix{
\tilde{Q}S(m)\ar[r]\ar[d]&\tilde{Q}R\ar[d]\\
\tilde{Q}D(m)\ar[r]\ar@{.>}[ur] &S
}
\end{equation}
To prove that $f$ is a weak equivalence, we must show that the following
map is an isomorphism for all $m$
\begin{equation}\label{fiso}
f:H^mN\tilde{Q}R\overset{\cong}\to H^mNS
\end{equation}
We first prove that \eqref{fiso} is surjective. If $x\in H^mNS$ is an element,
call $x$ the map $\Z\to H^mNS$, $1\mapsto x$. Choose a cochain homomorphism
$\hat{x}:\Z[m]\to NS$ inducing $x$. We have an exact sequence
\[
 \xymatrix{
0 \ar[r]& \Z[m+1]\ar[r]& \Z<m, m+1>\ar[r]&\Z[m] \ar[r]& 0}
\]
Because both $N$ and $Q$ are exact, we have a solid line commutative diagram 
\begin{equation}\label{quierache}
\xymatrix@1{
NQ \Z [m+1] \ar@{ >->}[d] \ar[rrr]^0  &                              & &N\tilde{Q}R \ar@{->>}[d]^{f}\\
NQ \Z <m, m+1> \ar[r] \ar@{.>}[urrr]^h & NQ \Z [m] \ar@{->>}[dr] \ar[rr]&                &NS  \\
                                    &                              & \Z [m] \ar[ur]^x &   
} 
\end{equation}
To prove that the dotted arrow exists, apply the functor $K$ to obtain a commutative
diagram
\begin{equation}\label{seteci}
\xymatrix{
Q\Z[m+1]\ar[r]^{Q(0)}\ar[d]&\tilde{Q}R\ar[d]\\
Q\Z<m,m+1>\ar[r]\ar@{.>}[ur]&S
}
\end{equation}
Next use Lemma \ref{qttq} to obtain a diagram of the form \eqref{liftexi}
in which the top row is in the image of $Q$, whence the dotted arrow
exists in \eqref{liftexi}, whence also in \eqref{seteci} and \eqref{quierache}.
Call $y$ the arrow $NQ\Z[m] \to NQR$ induced by $h$. 
Then the image of $1$ through $\bar y : \Z=H^m(NQ\Z[m])\to H^m N\tilde{Q}R$ maps to $x$
under \eqref{fiso}. This proves that \eqref{fiso} is surjective. To show it
is also injective, let $x:\Z[m]\to N\tilde{Q}R$ represent an element in the
kernel of \eqref{fiso}. Then $fx:\Z[m]\to NS$ factors through a map
$x':\Z<m-1,m>\to NS$.
Because $\pe$ is natural we have a commutative diagram 
\[
\xymatrix{
NQ \Z [m]\ar@{ >->}[d]\ar[r]^{x\pe} & N\tilde{Q}R \ar@{->>}[d]^{f}\\ 
NQ\Z<m-1, m> \ar[r]_(.7){x'\pe} &NS 
}
\] 
Because $\pe$ is an equivalence, it suffices to show that $x\pe$ induces
the zero map in cohomology. Next, by virtue of \ref{kequivq} there is a 
homotopy
$fx\pe\to fNQ(\overline{x\pe})$. Because $\quad$
\mbox{$Q\Z[m]\cofi Q\Z<m-1,m>$} 
is a cofibration this homotopy extends to one between $x'\pe$
and some map $y$ which fits into the following commutative diagram.
\begin{equation}\label{minino}
\xymatrix{
NQ \Z [m] \ar[r]^{NQ(\overline{x\pe})} \ar@{ >->}[d] &N\tilde{Q}R \ar@{->>}[d]^{f}\\
NQ\Z<m-1, m> \ar[r]_(.7){y} \ar@{.>}[ur]        &NS 
}
\end{equation}
The same argument used during the course of the proof of the surjectivity
of \eqref{fiso} shows that the dotted arrow exists. Hence $x\pe$ induces
the zero map in cohomology, since it is homotopic to $NQ(\overline{x\pe})$,
and the latter induces zero by \eqref{minino}.
\end{proof}

\begin{coro}
\label{main}
The functor $\mathbb{L}\tilde{Q}:\ho DGR^*\to \ho\ass^\Delta$ of 
{\rm \ref{qpreserve}} is
an equivalence of categories.
\end{coro}
\begin{proof}
Immediate from \ref{main1} and \cite{hov}, 1.3.13.
\end{proof}
\begin{coro}\label{comain}
Let $K:\cha\to \Ab^\Delta$ be the Dold-Kan functor. If $A\in DGR^*$, equip
$KA$ with the product \eqref{trinus}. Then the left derived functor
$\mathbb{L}K$ of $DGR^*\to \ass^\Delta$, $A\mapsto KA$ 
is a category equivalence $\ho DGR^*\weq\ho\ass^\Delta$.
\end{coro}
\begin{proof}
Immediate from \ref{main} and \ref{qpreserve} ii).
\end{proof}

\bigskip 

\noindent{\it Acknowledgement.}
We are indebted to the referee for several improvements over
the original manuscript, and especially for a simplification of the proofs
of Section \ref{equimain}. 
The second author wishes to acknowledge useful 
discussions with Juan Guccione and Gabriel Minian.

\end{document}